\documentclass[12pt]{article}
\usepackage{a4wide}

\usepackage{amsfonts}
\usepackage{amscd,amsmath}
\usepackage{euscript}

\usepackage[arrow,matrix,curve]{xy}\SilentMatrices
\def\xyma{\xymatrix@M.1em}

\usepackage{epsfig}
\usepackage{amssymb}
\newfont{\Ar}{msam10}
\newcommand{\MySquare}{\hfill {\mbox{\Ar\symbol{3}}}}
\usepackage{psfrag}
\newtheorem{theorem}{Theorem}

\newtheorem{lemma}[theorem]{Lemma}
\newtheorem{corollary}[theorem]{Corollary}
\newtheorem{definition}[theorem]{Definition}
\input xy
\xyoption{all}

\begin{document}

\title{A colimit of classifying spaces}

\author{Graham Ellis \\
{\small Mathematics Department, National University of Ireland,
Galway}
\medskip\\
Roman Mikhailov \\
{\small  Steklov Mathematical Institute, Moscow } }

\maketitle

\date{}

\abstract{ We recall a group-theoretic description of the first
non-vanishing homotopy group of a certain (n+1)-ad of spaces and
show how it yields several
  formulae for  homotopy and homology groups of specific spaces.
In particular we obtain an alternative proof of J. Wu's group-theoretic description of the homotopy groups of a $2$-sphere.
}

\begin{section}{Introduction}
 This article is motivated by  recent interest in describing
  homotopy and homology groups as intersections of subgroups (see for example
\cite{bak,bauesmikhailov,brown,casas,donadze,everaert,gutierrezratcliffe,wu}).
Some years ago the first author discovered a Hopf type formula for the
higher integral homology of a group \cite{ellisJPAA} using
F. Keune's
 theory of relative nonabelian derived functors \cite{keune}.
  He subsequently observed that the formula could 
also be obtained as an easy consequence of
a higher homotopy van Kampen theorem of R. Brown and J.-L.
Loday \cite{brownloday}
 and published this  topological proof
with R. Brown in \cite{brownellis}.
 In this paper we take a closer look at the techniques
underlying the topological
proof with a view to pursuing ideas on intersections in free groups
suggested
 by the second author and H.-J. Baues in \cite{bauesmikhailov}.

Let  $G$  be a group with normal subgroups $N_1,\ldots,N_n$.
 We are
interested in the $n$th and $(n+1)$st homotopy groups of the
topological space $X$ arising
 as the homotopy colimit of
classifying spaces
$B(G/\prod_{i\in I} N_i)$ where $I$ ranges over all
subsets $I \subsetneq \{1,\ldots,n\}$. We need a
 connectivity condition and define
 an $m$-tuple of normal subgroups $(N_1,\ldots,N_m)$  to be {\it connected} if either $m\le 2$ or
 $m\ge 3$ and for  all subsets $I, J \subset \{1,\cdots,m\}$ with $ |I|\ge 2,
|J|\ge 1$ the following equality holds:
\begin{equation}\left( \bigcap_{i\in I} N_i \right) \left( \prod_{j\in J}N_j\right) =
\bigcap_{i\in I} \left( N_i (\prod_{j\in J}N_j) \right)\,.\end{equation}

 Part of the statement of our main result involves
  a group  $T(N_1,\ldots,N_n)$
  which is recalled from \cite{ellissteiner}  
 in Theorem \ref{thmnad} below. Our main result is:
\begin{theorem}\label{mainthm} If the $(n-1)$-tuple $(N_1,\ldots,\hat N_i,\ldots,N_n)$ is connected for each $1\le i\le n$ then the
 above colimit  $X$ has
\begin{equation}\label{iso1}\pi_n(X) \cong \frac{N_1\cap \ldots \cap N_n}{\prod_{I\cup J=\{1,\ldots,n\},
 {I\cap J=\emptyset}}
[\cap_{i\in I}N_i,\cap_{j\in J}N_j]} , \end{equation}
\begin{equation}\label{iso2}
\pi_{n+1}(X) \cong \ker(\partial\colon T(N_1,\ldots,N_n)\rightarrow G) .
\end{equation}
\end{theorem}
 For $n=2$ isomorphism (\ref{iso1}) is given
in \cite{brown} and reads
\begin{equation}\pi_2(X) \cong \frac{N_1\cap N_2}{[N_1, N_2]}.\end{equation}
For $n=2$ isomorphism
(\ref{iso2})
is given in \cite{brownloday}. For $n=3$ isomorphism (\ref{iso1}) is seemingly
new and reads
$$\pi_3(X) = \frac{N_1\cap N_2\cap N_3}{[N_1,N_2\cap N_3][N_2,N_1\cap N_3][N_3,
N_1\cap N_2]}.$$
 We shall explain how
 isomorphism (\ref{iso1}) implies:  the group-theoretic description
of $\pi_n(S^2)$ given by J. Wu  \cite{wu};   homomorphisms
 of second and third homotopy groups
given in \cite{gutierrezratcliffe} and
  \cite{bauesmikhailov} respectively;
Hopf type formulae for the higher homology of a group
given in \cite{ellisJPAA, brownellis}; seemingly new
results on the existence of
torsion in certain groups arising from "almost aspherical"
presentations.

\medskip
Our approach is heavily influenced by work of R. Brown and J.-L. Loday.
In \cite{loday} Loday introduced a functor $\Pi\colon ($n-cubes\
of\ spaces$) \rightarrow ($cat$^n$-groups$)$ from topology to algebra.
 Brown and Loday
\cite{brownloday} proved that this functor preserves certain connectivity conditions and certain
colimits.  An explicit  description of colimits in the algebraic category
 can be found
 in \cite{brown} for $n=1$, in \cite{loday,brownloday} for $n=2$,
  and in
\cite{ellissteiner} for the general case $n\ge 1$.
 We shall
explain how Theorem \ref{mainthm} is a fairly immediate consequence of this body of work. Indeed, Theorem \ref{mainthm} is really just an observation, but as this observation seems not to have been made before, and as it has a number of useful consequences, we give it the status of a theorem.

  We begin with the
case of  three normal subgroups.
\end{section}

\begin{section}{The case $n=3$}
Suppose that $L, M, N $ are normal subgroups of a group $G$. This
data gives rise to a commutative cube of  spaces

$$
{ \xyma{& & B(G) \ar@{->}[rr] \ar@{->}[ldd]
\ar@{-}[dd] && B(G/L) \ar@{->}[ddd] \ar@{->}[ldd]\\ \\
& B(G/M) \ar@{->}[ddd] \ar@{->}[rr] & \ar@{->}[d] & B(G/LM) \ar@{->}[ddd]\\
& & B(G/N) \ar@{-}[r] \ar@{->}[ldd] & \ar@{->}[r] & B(G/LN)
\ar@{->}[ldd]\\ & & \\ & B(G/MN) \ar@{->}[rr] & & X  } }
$$
in which $B(G)$ denotes the classifying space of $G$ and $X$ is
the homotopy pushout of the diagram of classifying spaces. Working
up to homotopy type, we can extend this cube of spaces to a
diagram of 27 spaces
$$
{ \xyma{& & F_{-1,-1,-1} \ar@{->}[rr] \ar@{->}[ldd]
\ar@{-}[dd] && F_{0,-1,-1} \ar@{->}[ddd] \ar@{->}[ldd]\\ \\
& F_{-1,0,-1} \ar@{->}[ddd] \ar@{->}[rr] & \ar@{->}[d] & F_{0,0,-1} \ar@{->}[ddd]\\
& & F_{-1,-1,0} \ar@{-}[dd] \ar@{-}[r] \ar@{->}[ldd] & \ar@{->}[r]
& F_{0,-1,0} \ar@{->}[ldd]\\ & & \\ & F_{-1,0,0} \ar@{->}[rr] &
\ar@{->}[d] & B(G) \ar@{-}[dd] \ar@{->}[ddl]
\ar@{->}[rr] && B(G/L) \ar@{->}[ddd] \ar@{->}[ddl]\\ & & F_{-1,-1,1} \ar@{->}[ddl] \\
& &
B(G/M)\ar@{->}[ddd] \ar@{->}[rr] & \ar@{->}[d] & B(G/LM) \ar@{->}[ddd]\\
& F_{-1,0,1} \ar@{->}[ddl] \ar@{-}[r] & \ar@{->}[r] &  B(G/N)
\ar@{-}[r]
\ar@{->}[ddl] & \ar@{->}[r] & B(G/LN) \ar@{->}[ddl]\\
\\ F_{-1,1,1} \ar@{->}[rr] && B(G/MN) \ar@{->}[rr] && X } }
$$
in which each row and each column is a fibration sequence. For
precise details of the construction see \cite{loday,brownloday}.

\medskip
The homotopy exact sequence of the fibration $F_{-1,1,1}
\rightarrow B(G/MN) \rightarrow X$ gives isomorphisms
\begin{equation}\label{eq3} \pi_n(X)\cong \pi_{n-1}(F_{-1,1,1}) \ \ (n\ge 3).\end{equation}
The  exact sequence of the fibration $F_{-1,-1,1} \rightarrow
F_{-1,0,1} \rightarrow F_{-1,1,1}$ gives isomorphisms
\begin{equation}\pi_{2}(F_{-1,1,1}) \cong \ker (\pi_1(F_{-1,-1,1}) \rightarrow L/(L\cap N)) ,
\label{eq4}\end{equation}
\begin{equation}\pi_{n}(F_{-1,1,1}) \cong \pi_{n-1}(F_{-1,-1,1}) \ \ (n\ge 4).
\label{eq5}\end{equation}
Note that
\begin{equation}\pi_1(F_{-1,0,0}) \cong L\cap M .\label{eq6}\end{equation}

\noindent {\bf Notation.} For $x,y\in G$ we define
$[x,y]=xyx^{-1}y^{-1}$ and ${^x}y=xyx^{-1}$.
 We  write $A\sqcup B=<n>$ to mean that $A, B$ are nonempty subsets of
$\{1,\ldots,n\}$ with union equal to $\{1,\ldots,n\}$ and with
trivial intersection. We write $A\sqcup B\sqcup C =<n>$ to mean
that $A,B,C$ are nonempty sets with union equal to
$\{1,\ldots,n\}$ and with trivial pairwise intersections. We set
$N_1=L, N_2=M, N_3=N$ and $N_A=\cap_{i\in A}N_i$.

\medskip
 Brown and Loday's higher van Kampen theorem \cite{brownloday} asserts that
 a functor  $\Pi\colon ($3-cubes\
of\ spaces$) \rightarrow ($cat$^3$-groups$)$ preserves certain connectivity and colimits.
 In paticular it asserts that all
27 spaces in the diagram are path-connected.
  The following
consequence of the colimit property
 was proved in \cite{ellisthesis} and published in the form of a more general result
in \cite{ellissteiner}.

\begin{theorem}\cite{ellisthesis,ellissteiner}\label{thm1}  The group
$\pi_1(F_{-1,-1,-1})$ is isomorphic to the group $T(L,M,N)$
generated by symbols  $$a\otimes_{A,B} b$$ for all $A\sqcup B=
<3>,  a\in N_A, b\in  N_B$ ,
 subject to the  relations
$$a\otimes_{A,B} b = (b\otimes_{B,A} a)^{-1} ,$$
$$aa'\otimes_{A,B} b = ( {^a}a'\otimes_{A,B} {^a}b)(a\otimes_{A,B} b) ,$$
$$ ({^u}[u^{-1}, v]\otimes_{U\cup V,W} {^u}w)({^w}[w^{-1},u]\otimes_{W\cup U,V}
{^w}v)({^v}[v^{-1},w]\otimes_{V\cup W,U} {^v}u),$$
$$(a\otimes_{A,B} b)(a'\otimes_{A',B'} b')(a\otimes_{A,B}
b)^{-1}={^{[a,b]}}a'\otimes_{A',B'} {^{[a,b]}}b' $$ for  $A\sqcup
B= A'\sqcup B'=<3>,  a\in N_A, a'\in N_{A'}, b\in  N_B, b'\in
 N_{B'},
 U\sqcup V\sqcup W=<3>, u\in N_U, v\in N_V, w\in N_W.$
The homomorphism $T(L,M,N)\rightarrow \pi_1(F_{-1,-1,0})=L\cap M$
maps $x\otimes y$ to the commutator $[x,y]$; this homomorphism has
the structure of a crossed module
 $\partial\colon T(L,M,N)\rightarrow G, x\otimes y\mapsto [x,y]$
with action of $g\in G$ given by ${^g}(x\otimes y) = ({^g}x\otimes
{^g}y)$.
\end{theorem}

The homotopy exact sequence of the fibration $F_{-1,-1,-1}
\rightarrow F_{-1,-1,0} \rightarrow F_{-1,-1,1}$ together with
isomorphisms (\ref{eq3})-(\ref{eq6}) and Theorem \ref{thm1} yield the following
version of our main result for $n=3$.

\begin{theorem}\label{thm2}
There are isomorphisms
$$\pi_3(X) \cong \frac{L\cap M \cap N}{[L,M\cap N][M,L\cap N][N,L\cap M]},$$
$$\pi_4(X) \cong \ker(\partial\colon T(L,M,N)\rightarrow G) .$$
\end{theorem}

For completeness we also include
 isomorphisms for the first two homotopy groups of $X$; 
these  can both be deduced from Theorem \ref{thm1} using
 the path-connectivity
of the 27 spaces and the exact sequence of a 
fibration.
\begin{equation}\pi_1(X)\cong \frac{G}{LMN},\end{equation}
\begin{equation}\pi_2(X) \cong \frac{LM\cap MN}{M(L\cap N)}.\label{pitwo}\end{equation}
Note that the left-hand side of (\ref{pitwo}) is symmetric in $L,M,N$. Hence the right-hand side is also symmetric in these three arbitrary normal subgroups of $G$. 

The first isomorphism of the following corollary of Theorem \ref{thm2}
 is a result of Wu
\cite{wu}.

\begin{corollary}\label{cor3}
Let $G=H\ast H$ be the free product of two copies of a group $H$.
Let $\lambda, \mu, \nu \colon H\ast H\rightarrow H$ be the first
projection, second projection and multiplication homomorphisms
respectively. Set $L=\ker \lambda$, $M=\ker \mu$, $N=\ker \nu$.
Then the suspension $SK(H,1)$ of a classifying space for $H$ has
$$\pi_3(SK(H,1)) = \frac{L\cap M \cap N}{[L,M\cap N][M,L\cap N][N,L\cap M]},$$
$$\pi_4(SK(H,1)) \cong \ker(\partial\colon T(L,M,N)\rightarrow G) .$$

\end{corollary}

\noindent {\bf Proof.} Let $B$ be a classifying space for $H$.
Then the wedge $B\vee B$ is a classifying space of $G=H\ast H$.
Let $CB$ denote the cone on $B$. Then $CB\vee B$ is a classifying
space for $G/L$, $B\vee CB$ is a classifying space for $G/M$, and
$ B\times [0,1]$ is a classifying space for $G/M$. Moreover, the
pushout (or union) of these classifying spaces is a suspension
$X=SK(H,1)$. The corollary follows from Theorem \ref{thm2}.
\MySquare

\bigskip
The first homomorphism of the following corollary was proved by
Baues and Mikhailov \cite{bauesmikhailov}.

\begin{corollary}\label{cor4}
Let $K_1,K_2,K_3$ be CW-subspaces of a connected CW-space $K$ for
which $K=K_1\cup K_2\cup K_3$ and $K^1=K_1\cap K_2\cap K_3$ is the
1-skeleton of $K$. Let $F=\pi_1K^1$ and $R_i=\ker (\pi_1K^1
\rightarrow \pi_1K_i)$ ($i=1,2,3$).
 Then there are homomorphisms
$$\alpha_3\colon \pi_3(K) \longrightarrow \frac{R_1\cap R_2\cap R_3}{[R_1,R_2\cap
R_3][R_2,R_1\cap R_3][R_3,R_1\cap R_2]} ,$$
$$\alpha_4\colon \pi_4(K) \longrightarrow \ker(\partial\colon
T(R_1,R_2,R_3)\rightarrow F) .$$ In some cases both $\alpha$ and
$\beta$ are isomorphisms.
\end{corollary}

\noindent{\bf Proof.} Each homomorphism $\pi_1K^1 \rightarrow
\pi_1K_i$ is surjective since the 1-skeleton of $K_i$ equals the
1-skeleton $K^1$ of $K$. Hence $\pi_1K_i\cong F/R_i$. We can
identify $K^1$ with the classifying space $B(F)$. There are maps
$\iota_i\colon K_i\rightarrow B(F/R_i)$ inducing  isomorphisms of
fundamental groups. The maps $\iota_i$ induce a map $\iota\colon
K=K_1\cup K_2\cup K_3 \rightarrow X$ where $X$ is the colimit of
Theorem \ref{thm1} with $G=F, L=R_1, M=R_2, N=R_3$. The required
homomorphisms $\alpha_k$ are obtained by composing the
isomorphisms of Theorem \ref{thm1} with the induced homomorphisms
$\pi_k(\iota)\colon \pi_k(K)\rightarrow \pi_k(X)$ ($k=3,4$).

Consider the 2-sphere $K=SK(\mathbb Z,1)$ expressed as a union of
$K_1=CS^1\vee S^1$, $K_2=S^1\vee CS^1$, $K_3=S*1\times [0,1]$ as
in the proof of Corollary \ref{cor3}.  In this case the maps
$\alpha_3, \alpha_4$ are isomorphisms.
 \MySquare

\bigskip

One can obtain group-theoretic results using Theorem \ref{thm2}.
Let $\cal X$ be an alphabet and $F=F(X)$ be the free group
generated by $\cal X$. Let $r_1$,\ldots ,$r_n$ be words in $F$.
The
 group presentation
${\cal P}(X,\ r_1,\ldots, r_n)=\langle X\ |\ r_1,\ldots
,r_n\rangle$ is said to be {\it aspherical} if the associated
2-complex has trivial second homotopy group. We call the
presentation {\it almost aspherical} if for every proper subset
$\cal S$ of the set $\{r_1$,\ldots ,$r_n\}$ the presentation
$\langle X | \cal S\rangle$ is aspherical.
 For example, the 4-string braid group has the almost aspherical presentation
$\langle x,y,z | xyx=yxy, yzy=zyz, xz=zx\rangle$.

\begin{corollary}\label{cor4A}
Let ${\cal P}(X,r_1, r_2, r_3)$ be an almost aspherical
presentation. Let $R_i$ be the normal closure in $F=F(X)$ of the
relator ${\cal R}_i$. Then the quotient group
$$\frac{F}{[R_1,R_2\cap R_3][R_2,R_1\cap R_3][R_3,R_1\cap R_2]}$$
is torsion free.
\end{corollary}

\noindent{\bf Proof.} The 2-complexes associated to the
presentations $\langle X\ |\ r_i\rangle$ and $\langle X\ |\ r_i,\
r_j\rangle$ are  classifying spaces for $F/R_i$ and $F/R_iR_j$
($1\le i\ne j \le 3$). So on taking $G=F, L=R_1, M=R_2, N=R_3$ in
Theorem \ref{thm2} we can identify the colimit space $X$ with the
2-complex associated to the three-relator presentation. Theorem
\ref{thm2} provides an exact sequence
$$0\rightarrow \pi_3(X) \rightarrow \frac{F}{[R_1,R_2\cap R_3][R_2,R_1\cap
R_3][R_3,R_1\cap R_2]} \rightarrow \frac{F}{R_1\cap R_2\cap R_3}
\rightarrow 1.$$
 The universal cover of $X$ is homotopy equivalent to a wedge of 2-spheres and thus
$\pi_3(X)$ is known to be torsion free. The group $F/R_1\cap
R_2\cap R_3$ maps into the  group $(F/R_1)\times (F/R_2) \times
(F/R_3)$ and each $F/R_i$ is torsion free since it has an
aspherical presentation. The corollary follows. \MySquare

\bigskip

It is easily shown that $T(G,G,G)$ is isomorphic to the nonabelian
symmetric square $G\tilde \otimes G$ introduced by  Dennis
\cite{dennis} and then Brown and Loday \cite{brownloday}. In light
of this isomorophism, the following corollary can be found in
\cite{brownloday}.

\begin{corollary}\label{cor5}
For the double suspension $S^2K(G,1)$ of a classifying space
$K(G,1)$ there are isomorphisms
$$\pi_3(S^2K(G,1)) \cong \frac{G}{[G,G]} ,$$
$$\pi_4(S^2K(G,1)) \cong \ker(\partial\colon T(G,G,G)\rightarrow G) .$$
\end{corollary}

\noindent{\bf Proof.} It suffices to observe that, in the case
$L=M=N=G$,
 the colimit $X$ is of the homotopy type of $S^2K(G,1)$.
\MySquare

\bigskip
The above results on homotopy groups rely heavily on the homotopy
exact sequence of a fibration. For calculations in group
cohomology one defines $B(G,N)$ to be the homotopy cofibre of the
map $B(G)\rightarrow B(G/N)$ and define $H_k(G)=H_k(B(G),\mathbb
Z)$, $H_k(G,N)=H_{k+1}(B(G,N),\mathbb Z)$. This yields a long
exact sequence
$$\cdots \rightarrow H_{k+1}(G/N) \rightarrow H_{k+1}(G,N) \rightarrow H_k(G)
\rightarrow H_k(G/N) \rightarrow \cdots .$$
 One defines $B(G,M,N)$ to be the homotopy cofibre of the canonical map
$B(G,N) \rightarrow B(G/M,NM/M)$ and $H_k(G,M,N) =
H_{k+2}(B(G,M,N),\mathbb Z)$. There is thus a long exact sequence
$$\cdots \rightarrow H_{k+1}(G/M,MN/M) \rightarrow H_{k+1}(G,M,N) \rightarrow
H_k(G,N) \rightarrow H_k(G/M,MN/N) \rightarrow \cdots . $$

\medskip
The  isomorphism in the following corollary was originally proved
by a purely algebraic argument using the theory of nonabelian left
derived functors \cite{ellisJPAA}; however, the first proof to
appear in print was a short
 topological one \cite{brownellis}. We recall the topological proof here.

\begin{corollary}\label{cor6}
Let $R$ and $S$ be normal subgroups of a free group $F$ and
suppose that both quotients $F/R$ and $F/S$ are free. Set
$G=F/RS$. Then
$$H_3(G)=\frac{R\cap S\cap [F,F]}{[R,S][R\cap S,F]} .$$
\end{corollary}

\noindent{\bf Proof.} It is readily seen that, for any normal
subgroups $M,N$ in $G$,
 the
cofibre $B(G,M,N)$ is homotopy equivalent to the above colimit $X$
in the case $L=G$. Theorem \ref{thm2} shows that
$$\pi_3(B(G,M,N)) \cong  \frac{M\cap N}{[G,M\cap N][M,N]}.$$
Since $\pi_k(B(G,M,N))=0$ for $k=1,2$, the Hurewicz isomorphism
$\pi_3(B(G,M,N)) \cong H_3(B(G,M,N),\mathbb Z)$ gives
$$H_1(G,M,N) \cong  \frac{M\cap N}{[G,M\cap N][M,N]}.$$
Using $H_k(F)=H_k(F/R)=H_k(F/S)=0$ for $k\ge 2$, the  exact
homology sequences of a cofibration yield an isomorphism
$$H_3(G)\cong \ker(H_1(F,R,S)\rightarrow H_1(F)) .$$
This yields the corollary. \MySquare

\bigskip
The above formula for $H_1(G,M,N)$ is of independent interest. For
a group-theoretic description of $H_2(G,M,N)$ we define the group
$E(L,M,N)$ to be the quotient of $T(L,M,N)$ obtained by imposing
the extra relations
$$x\otimes x = 1$$
for all $x\in L\cap M\cap N$. There is an induced
 crossed module $\partial\colon E(L,M,N)\rightarrow G$. The following result was
proved algebraically in \cite{ellisJPAA}. To keep things
topological we outline a topological proof of the second
isomorphism.

\begin{corollary}\label{cor7}
For a group $G$ with normal subgroups $M,N$ there are isomorphisms
$$H_1(G,M,N) \cong \frac{M\cap N}{[G,M\cap N][M,N]} ,$$
$$H_2(G,M,N) \cong \ker(E(G,M,N) \rightarrow G) .$$
\end{corollary}

\noindent{\bf Proof.} Let $L=G$. We have already derived the first
isomorphism from Theorem \ref{thm2}. For the second we consider
J.H.C. Whitehead's exact sequence
$$H_5(X) \rightarrow \Gamma_4(X) \stackrel{\eta}{\rightarrow} \pi_4(X) \rightarrow
H_4(X) \rightarrow 0 .$$
 A description of $\Gamma_4(X)$ is given in \cite{baues}; in our case we find
$$\Gamma_4(X) = ((M\cap N)/[G,M\cap N][M,N]) )\otimes \mathbb
Z_2.$$
 To complete the proof one has to check that $\eta$ is induced by the
homomorphism $M\cap N \rightarrow T(G,M,N), x\mapsto x\otimes x$.
\MySquare
\end{section}

\begin{section}{The case $n\ge 2$}
We need to consider the category $\{0<1\}$ with two objects and
one non-identity arrow $0\rightarrow 1$. Following \cite{loday} we
define an {\it $n$-cube of spaces} to be a functor $ S\colon
\{0<1\}^n\rightarrow ($Spaces$)$. We also need to consider
  the category $\{-1<0<1\}$ with three objects and two non-identity arrows
$-1\rightarrow 0$ and $0\rightarrow 1$. An object in this category
is an n-tuple $\Delta\in \{0,1\}^n$; we define a {\it row} to
consist of a sequence of arrows $\Delta'\rightarrow
\Delta\rightarrow \Delta"$ where the $n$-tuples $\Delta', \Delta,
\Delta"$ are identical except in say the $i$th coordinate.
 Following \cite{loday} we define
 an {\it $n$-cube of fibrations} to be a functor
$F\colon \{-1<0<1\}^n\rightarrow ($Spaces$)$ that maps each row to
a fibration sequence.
 It is explained in \cite{loday,brownloday} how an $n$-cube of spaces $S$ gives rise to an
$n$-cube of fibrations $F= \overline S$ where, up to homotopy type
of maps, the restriction of $\overline S$ to $\{0<1\}^n$ is $S$.

\begin{definition}\label{connect}
 Following \cite{brownloday} we say that an $n$-cube of spaces $S$
is connected if each space in the $n$-cube of fibrations
$\overline S$ is path-connected.
\end{definition}

Given an $n$-tuple $\Delta\in \{0,1\}^n$ we shall write $i\in
\Delta$ to mean that the $i$th coordinate of $\Delta$ eqials $1$.
 To $n$ normal subgroups $N_1,\ldots,N_n$ of a group $G$ we associate
 the $n$-cube of
spaces $B=B(G,N_1,\ldots,N_n)\colon \{0<1\}^n\rightarrow
($Spaces$)$ in which
$$B_\Delta = B(\frac{G}{\prod_{i\in \Delta}N_i})$$  is the classifying space of
 the quotient group $G/\prod_{i\in \Delta}N_i$,
and the maps are induced by the canonical inclusions.

The verification of the following lemma, which uses
 the homotopy exact sequence of a fibration,  is left to the reader.

\begin{lemma}\label{lem10} Let $B=B(G,N_1,\ldots,N_n)$ be the
 $n$-cube of spaces  associated to $n$ normal subgroups
 of $G$.\newline
\begin{enumerate}
\item[(i)] The n-cube $B$ is connected if $n=1,2$. \item[(ii)] For
$n\ge 3$ the $n$-cube $B$
 is  connected  if and only if
for all subsets $I, J \subset \{1,\cdots,n\}$ with $ |I|\ge 2,
|J|\ge 1$ the following equality holds:
$$\left( \bigcap_{i\in I} N_i \right) \left( \prod_{j\in J}N_j\right) =
\bigcap_{i\in I} \left( N_i (\prod_{j\in J}N_j) \right)\,.$$
\end{enumerate}
\end{lemma}

Let $B=B(G,N_1,\ldots,N_n)$. Let $\Delta_{\rm final}
=(1,1,\ldots,1)$ denote the $n$-tuple with each coordinate equal
to $1$.
 Let $\Delta_{\rm initial} =(-1,-1,\ldots,-1)$
denote the $n$-tuple with each coordinate equal to $-1$.
 We are interested in the space $X$ arising as the colimit of classifying
spaces $\{B_\Delta : \Delta\in \{0<1\}^n, \Delta\ne\Delta_{\rm
final}\}$.
 Define $C=C(G,N_1,\ldots,N_n)\colon\{0<1\}^n\rightarrow ($Spaces$)$
 to be the $n$-cube of spaces obtained from $B$ by replacing the
classifying space $B_{\Delta_{\rm final}}$ with the colimit $X$.
Let $F=\overline C$ be the $n$-cube of fibrations associated to
$C$. The following theorem is a  special case of Theorem 2.3 in
 \cite{ellissteiner} and is a consequence of Brown and Loday's higher van
Kampen theorem.

\begin{theorem}\cite{ellissteiner}\label{thmnad} Let $N_1,\ldots,N_n$ be
 normal subgroups of $G$ such that the (n-1)-cube of
spaces $B(G,N_1,\ldots,\hat N_i,\dots,N_n)$ is connected for each
$1\le i\le n$. Let $F=\overline C$ be the $n$-cube of fibrations
associated to $C=C(G,N_1,\ldots,N_n)$. Every space in $F$ is
path-connected and
  $\pi_1(F_{\Delta_{\rm initial}})$ is isomorphic to the group $T(N_1,\ldots,N_n)$
generated by symbols
 $$a\otimes_{A,B} b$$
for all $A\sqcup B= <n>,  a\in N_A, b\in  N_B$ ,
 subject to the  relations
$$a\otimes_{A,B} b = (b\otimes_{B,A} a)^{-1} ,$$
$$aa'\otimes_{A,B} b = ( {^a}a'\otimes_{A,B} {^a}b)(a\otimes_{A,B} b) ,$$
$$ ({^u}[u^{-1}, v]\otimes_{U\cup V,W} {^u}w)({^w}[w^{-1},u]\otimes_{W\cup U,V}
{^w}v)({^v}[v^{-1},w]\otimes_{V\cup W,U} {^v}u),$$
$$(a\otimes_{A,B} b)(a'\otimes_{A',B'} b')(a\otimes_{A,B}
b)^{-1}={^{[a,b]}}a'\otimes_{A',B'} {^{[a,b]}}b' $$ for  $A\sqcup
B= A'\sqcup B'=<n>,  a\in N_A, a'\in N_{A'}, b\in  N_B, b'\in
 N_{B'},
 U\sqcup V\sqcup W=<n>, u\in N_U, v\in N_V, w\in N_W.$
The homomorphism $T(N_1,\ldots,N_n)\rightarrow
\pi_1(F_{-1,-1,\ldots,-1,0})\le G$ maps $x\otimes y$ to the
commutator $[x,y]$; this homomorphism has the structure of a
crossed module
 $\partial\colon T(N_1,\ldots,N_n)\rightarrow G, x\otimes y\mapsto [x,y]$
with action of $g\in G$ given by ${^g}(x\otimes y) = ({^g}x\otimes
{^g}y)$.
\end{theorem}

\noindent{\bf Proof of Theorem \ref{mainthm}}.
Our derivation of  Theorem \ref{thm2} from Theorem \ref{thm1}
extends routinely to yield our main Theorem \ref{mainthm}
as a  corollary of Theorem \ref{thmnad}. One simply has to apply the homotopy exact sequence of a fibration $n$ times. Lemma \ref{lem10} provides the algebraic version of the connectivity condition. \MySquare

\bigskip
 As already mentioned, Theorem \ref{mainthm} for  the case $n=2$ is given in
\cite{brownloday}. This case yields the following earlier result
of Brown \cite{brown}

\begin{corollary}\cite{brown}\label{cor13} Given two normal subgroups $M, N$ in a
group $G$ we have
$$\pi_2(B(G/M) \cup_{B(G)} B(G/N)) = \frac{M\cap N}{[M,N]}.$$
\end{corollary}

 Corollary \ref{cor13} implies that
 $R_i\cap R_j=[R_i,R_j]\ (i\ne j)$ in
Corollary \ref{cor4A}. Corollary \ref{cor13} also imples that
$L\cap M=[L,M], M\cap N=[M,N], L\cap N=[L,N]$ in Corollary
\ref{cor3}. Furthermore, Corollary \ref{cor13} can be used to recover the exact sequence of second homotopy groups proved by M. Gutierrez and J. Ratcliffe in
\cite{gutierrezratcliffe}.

\medskip
With one exception, the corollaries in Section 3 and their proofs
extend to analogous results involving $n\ge 2$ normal subgroups.
The exception is Corollary 8 where we we do not know an explicit
description of Whitehead's sequence in higher dimensions necessary
for extending the topological proof.

\end{section}

\begin{section}{Group-theoretical applications}
Some unexpected results in group theory  arise from
well-known results in homotopy theory. 
 For example Corallary \ref{cor4A}
generalizes to the following.

\begin{corollary}\label{cor14}
Let ${\cal P}(X,r_1,\ldots,
 r_n)$ be an almost aspherical presentation.
Let $R_i$ be the normal closure in $F=F(X)$ of the relator $R_i$.
Then the quotient group
$$\frac{F}{
\prod_{I\cup J=\{1,\ldots,n\},\ I\cap J=\emptyset} [\cap_{i\in
I}R_i,\cap_{j\in J}R_j] }$$ can only have $p$-torsion if the $n$th
homotopy group of a wedge of spheres can have $p$-torsion. So for
example, by a result of Serre, the group can have no $p$ torsion
for primes $p> 2n$.
\end{corollary}
{\bf Proof \ \ }
Consider a proper subset $I\subsetneq \{1,\ldots,n\}$. Since the group
presentation ${\cal P}(X,r_i \ (i\in I))$ is aspherical the associated 
$2$-complex $ K_I$ is a classifying space for the group $F/\prod_{i\in I} R_i$. Also, this space $K_I$ is the homotopy colimit (i.e. union) of the 
diagram of spaces
$\{ B(F/\prod_{i\in J}R_j) : J\subsetneq I\}=
\{ K_J : J\subsetneq I\}
$. Using induction on the size of $I$ together with
 the connectivity assertion of Theorem \ref{thmnad} we obtain that
the $|I|$-cube of spaces  $B(F,R_i \ (i\in I))$ is connected.  In particular each $(n-1)$-cube  $B(F,R_1,\ldots,\hat R_i, \ldots, R_n)$ is connected and so we can apply Theorem \ref{mainthm} to the normal subgroups $(F,R_1,\ldots,R_n)$;   
 we get a formula for the n-th homotopy group of a colimit space which, in this case, is the $2$-complex $K$ associated to the presentation
 ${\cal P}(X,r_1,\ldots,
 r_n)$. Now
 $\pi_2K \cong \pi_2\tilde K$ and $\tilde K$ is a wedge of $2$-spheres. 
 The remainder of the proof is now just a copy of the proof of Corollary 
 \ref{cor4A}.
\MySquare

\medskip
 For a second example we recall   
 a description of the homotopy groups of the 2-sphere due to
Wu \cite{wu}. 
 Let $F[S^1]$ be Milnor's $F[K]$-construction applied
to the simplicial circle $S^1$. This is the free simplicial group
with $F[S^1]_n$ a free group of rank $n\geq 1$ with generators
$x_0,\dots, x_{n-1}$. Changing the basis of $F[S^1]_n$ in the
following way: $y_i=x_ix_{i+1}^{-1},\ y_{n-1}=x_{n-1}$, we get
another basis $\{y_0,\dots,y_{n-1}\}$ in which the simplicial maps
can be written more easily. A combinatorial group-theoretical argument
then gives a description of the $n$-th homotopy group
of the loop space $\Omega\Sigma S^1$, which is isomorphic to the
homotopy group of $\pi_{n+1}(S^2)$ (see \cite{wu} for precise
details 
); one finds Wu's formula 
\begin{equation}\label{wujie}
\pi_{n+1}(S^2)\cong \frac{\langle y_{-1}\rangle^F\cap \langle
y_0\rangle^F \cap\dots\cap \langle
y_{n-1}\rangle^F}{[[y_{-1},y_0,\dots, y_{n-1}]]},\ n\geq 1
\end{equation}
where $F$ is a free group with generators $y_0,\dots, y_{n-1},$
$y_{-1}=(y_0\dots y_{n-1})^{-1}$, and  $[[y_{-1},y_0,\dots,
y_{n-1}]]$ is the normal closure in $F$ of the set of left-ordered
commutators \begin{equation}\label{comm}
[z_1^{\varepsilon_1},\dots, z_{t}^{\varepsilon_t}]\end{equation}
with the properties that $\varepsilon_i=\pm 1$,
$z_i\in\{y_{-1},\dots, y_{n-1}\}$ and all elements in
$\{y_{-1},\dots, y_{n-1}\}$ appear at least once in the sequence
of elements $z_i$ in (\ref{comm}).

Consider $n\geq 1$ and Milnor's construction $F[S^n]$. The lower
terms of the simplicial group $F[S^n]$ are the following:
\begin{align*}
& F[S^n]_n=F(\sigma),\\
& F[S^n]_{n+1}=F(s_0\sigma,\dots, s_n\sigma),\\
& F[S^n]_{n+2}=F(s_js_i\sigma\ |\ i<j),\\
& \dots
\end{align*}
For $n\geq 1,$ the generator of
$$\pi_{n+1}(F[S^n])\simeq \pi_{n+1}(\Omega(S^{n+1}))\simeq \pi_{n+2}(S^{n+1})$$
can be chosen as the commutator $[s_0\sigma,s_1\sigma]$ in
$F[S^n]_{n+1}$ (see \cite{wu}). Now these commutators considered
as elements of $F[S^n]_{n+1}$ define the elements from $F[S^1]_k$,
which correspond the homotopy classes of composition maps
$$
S^{k+1}\to S^k\to \dots\to S^3\to S^2,
$$
where every map is viewed as a suspension over the Hopf fibration.
Let us consider these elements.

\noindent 1. First, let $F_2=F(y_0,y_1)$, then the element
$$
[y_0,y_1]\notin [[y_{-1},y_0,y_1]]
$$
corresponds to the homotopy class of the Hopf fibration $S^3\to
S^2$.\\ \\
2. Let $F_3=F(y_0,y_1,y_2)$, then the element
$$
[[y_0,y_1],[y_0,y_1y_2]]\notin [[y_{-1},y_0,y_1,y_2]]
$$
corresponds to the homotopy class of the composition map $S^4\to
S^3\to S^2$.\\ \\
3. Let $F_4=F(y_0,y_1,y_2,y_3)$, then the element
$$
[[[y_0,y_1],[y_0,y_1y_2]],[[y_0,y_1],[y_0,y_1y_2y_3]]]\notin
[[y_{-1},y_0,y_1,y_2,y_3]]
$$
corresponds to the homotopy class of the composition map $S^5\to
S^4\to S^3\to S^2.$\\ \\
4. Let $F_5=F(y_0,y_1,y_2,y_3,y_4)$, then the element
\begin{multline*}
[[[[y_0,y_1],[y_0,y_1y_2]],[[y_0,y_1],[y_0,y_1y_2y_3]]],\\
[[[y_0,y_1],[y_0,y_1y_2]],[[y_0,y_1],[y_0,y_1y_2y_3y_4]]]]\\
\notin [[y_{-1},y_0,y_1,y_2,y_3,y_4]]
\end{multline*}
corresponds to the homotopy class of the composition map $S^6\to
S^5\to S^4\to S^3\to S^2$.\\ \\
5. Let $F_6=F(y_0,y_1,y_2,y_3,y_4,y_5)$, then the element
\begin{multline*}
[[[[[y_0,y_1],[y_0,y_1y_2]],[[y_0,y_1],[y_0,y_1y_2y_3]]],\\
[[[y_0,y_1],[y_0,y_1y_2]],[[y_0,y_1],[y_0,y_1y_2y_3y_4]]]],\\
[[[[y_0,y_1],[y_0,y_1y_2]],[[y_0,y_1],[y_0,y_1y_2y_3]]],\\
[[[y_0,y_1],[y_0,y_1y_2]],[[y_0,y_1],[y_0,y_1y_2y_3y_4y_5]]]]]\\
\in [[y_{-1},y_0,y_1,y_2,y_3,y_4,y_5]]
\end{multline*}
corresponds to the trivial homotopy class of the composition map
\begin{equation}\label{7sphere} S^7\to S^6\to S^5\to S^4\to S^3\to S^2.\end{equation} The
triviality of this map can be proved using standard methods in
homotopy theory \cite{Toda:62}. This is the simplest case of the
Nilpotence Theorem due to Nishida \cite{nishida}, which states
that every element in the ring of stable homotopy groups of
spheres is nilpotent.

\bigskip
We now turn to a group-theoretic application of Theorem \ref{mainthm}.
Consider a free group $F(X)$ with a basis $X$ and let $\mathcal
S=\{\mathcal R_1,\dots, \mathcal R_n\}$ be a collection of sets $\mathcal R_i$
 of
words, such that:
1) the presentation \begin{equation}\label{saa}\langle X\ |\
\mathcal R_1,\dots, \mathcal
R_n\rangle\end{equation} is not aspherical,
2) the presentation \begin{equation}\label{subpre}\langle X\ |\
\mathcal R_{i_1},\dots, \mathcal R_{i_k}\rangle\end{equation} is
aspherical for every proper subset $\{\mathcal R_{i_1},\dots,
\mathcal R_{i_k}\}$ of $\mathcal S$. Such a presentation we shall
call {\it $\mathcal S$-almost aspherical}. Observe that every
almost aspherical presentation is $\mathcal S$-almost aspherical,
where $\mathcal S$ is the collection of all relators, i.e. the
sets $\mathcal R_i$ contain single relators for all $i$.

Denote by $R_i$ the normal closure of $\mathcal R_i$ for
$i=1,\dots, n$. Observe that the standard 2-complex, associated to
an $\mathcal S$-almost aspherical presentation (\ref{saa}) is
homotopically equivalent to the homotopy colimit $X$ of the
$n$-cube of standard 2-complexes associated with subpresentations
of the type (\ref{subpre}), i.e. to the colimit of the non-final spaces
in the $n$-cube
$B(F,R_1,\dots,R_n)$. The asphericity condition 2) implies that
the above $n$-cube is $n$-connected and hence the $n$-th homotopy
group $\pi_n(X)$ can be computed by formula given in Theorem
\ref{mainthm}. However, the colimit $X$ is homotopically
equivalent to a (non-empty) wedge of 2-spheres and, therefore,
 formulas for the homotopy groups of  wedges of
2-spheres follow.

As an example, consider an almost aspherical presentation of the
trivial group:
$$
\mathcal P=\langle x_1,\dots, x_n\ |\ r_1,\dots, r_{n+1}\rangle
$$
with $n$ generators and $n+1$ relators $(n\geq 1)$. Let
$K_{\mathcal P}$ be the standard 2-complex associated with
$\mathcal P$. The chain complex $C_*K_{\mathcal P}$ has
$C_2K_{\mathcal P}=\mathbb Z^{n+1}$, $C_1K_{\mathcal P}=\mathbb
Z^n$ and natural epimorphism $C_2K_{\mathcal P}\to C_1K_{\mathcal
P}$, which kernel is isomorphic to $H_2(K_{\mathcal
P})=\pi_2(K_{\mathcal P})=\mathbb Z$, hence $K_{\mathcal P}$ is
homotopically equivalent to the 2-sphere $S^2$. Theorem
\ref{mainthm} implies the following presentation of the $(n+1)$-st
homotopy group of $S^2$:
\begin{equation}\label{2sphere}
\pi_{n+1}(S^2)\simeq \frac{R_1\cap \dots\cap R_{n+1}}{\prod_{I\cup
J=\{1,\ldots,n\},\ I\cap J=\emptyset} [\cap_{i\in I}R_i,\cap_{j\in
J}R_j]},
\end{equation}
where $R_i=\langle r_i\rangle^F,$ $F$ being the free group with
generators $x_1,\dots, x_n$. Clearly, in this case, the subgroup
of $F/({\prod_{I\cup J=\{1,\ldots,n\},\ I\cap J=\emptyset}
[\cap_{i\in I}R_i,\cap_{j\in J}R_j]})$ given in the right hand
side of (\ref{2sphere}) is central. In particular, consider the
presentation
$$
\mathcal P_{WU}=\langle x_1,\dots, x_n\ |\ x_1,\dots, x_n,
x_1\cdots x_n\rangle.
$$
It is easy to see that $\mathcal P_{WU}$ is almost aspherical. A
standard commutator calculus argument, given essentially in
Corollary 3.5 of \cite{wu} shows that
\begin{equation}
[[x_1\dots x_n,x_1,\dots, x_n]]=\prod_{I\cup J=\{1,\dots, n+1\},\
I\cap J=\emptyset}[\bigcap_{i\in I}R_i,\bigcap_{j\in J}R_j],
\label{comm2} \end{equation}
and Wu's isomorphism (\ref{wujie}) follows from (\ref{2sphere}). (One could also obtain (\ref{comm2}) using Theorem \ref{mainthm}.)

Now consider an arbitrary set of elements $r_1,\dots,r_{n+1}$ in
the free group $F$ with basis $x_1,\dots, x_n$ and with the
following property: the groups defined by presentations with $n$
generators and $n$ relators
$$
\langle x_1,\dots, x_n\ |\ r_1,\dots, r_{j-1},\hat
r_j,r_{j+1},\dots, r_{n+1} \rangle
$$
are trivial for all $j=1,\dots,n+1$ (i.e. we consider the
presentation without symbol $r_j$). Since every balanced
presentation, i.e. a presentation with equal number of generators
and relators, of the trivial group is aspherical, we have the
following: {\it either $\pi_{n+1}(S^2)$ can be presented as
(\ref{2sphere}) for $R_i=\langle r_i\rangle^F,\ i=0,\dots, n+1$,
or there is a counter-example to Whitehead asphericity conjecture,
i.e. there exists an aspherical 2-dimensional complex with
non-aspherical subcomplex.}

Since every subcomplex of an aspherical complex, with a single
2-cell, is aspherical (see \cite{A}), we obtain the following. Let
$F$ be a free group with generators $x_1,x_2$ and let
$r_1,r_2,r_3$ be words in $F$ such that the groups $F/R_1R_2,
F/R_1R_3, F/R_2R_3$ are trivial for $R_i=\langle r_i\rangle^F,\
i=1,2,3$. Then there is an isomorphism
$$
\frac{R_1\cap R_2\cap R_3}{[R_1,R_2\cap R_3][R_2,R_1\cap
R_3][R_3,R_1\cap R_2]}\simeq \mathbb Z.
$$
As an example, consider
$$r_1=x_1^2x_2^{-3},\ r_2=x_1^3x_2^{-4},\ r_3=x_1x_2x_1x_2^{-1}x_1^{-1}x_2^{-1}.$$
The groups
$$
\langle x_1,x_2\ |\ x_1^nx_2^{-(n+1)},\ x_1x_2x_1x_2^{-1}x_1^{-1}x_2^{-1}\rangle\\
$$ are trivial for $n\geq 2$. The above presentation of trivial groups is due to
Akbulut-Kirby \cite{ac}. It is easy to see that the presentation
$$
\langle x_1,x_2\ |\ x_1^2x_2^{-3},\ x_1^3x_2^{-4} \rangle
$$
defines the trivial group. It would be interesting to find a
natural generalization of this presentation to the case of $n$
generators and $n+1$ relators and to consider group-theoretical
meaning of the fact that the composition map (\ref{7sphere}) is
homotopically trivial.

\end{section}

\begin{section}{Other varieties}
The   Hopf type formula of Corollary \ref{cor6} and its higher
dimensional analogues were  first proved using a modification of a
technique of F. Keune for nonabelian derived functors. The
technique is very general and can be applied in situations such as
the homology of Lie algebras or varietal Baer invariants of
groups. Indeed, it was used in \cite{burnsellis} to derive a Hopf
type formula for  the  Baer invariant correcponding to the variety
of $2$-nilpotent groups; this formula was implemented on a
computer
 and applied to all groups of order at most 30.

The nonabelian
  derived functor approach and the generalised van Kampen theorem approach to Hopf
type formulae both apply only to connected $n$-cubes (see
Definition \ref{connect} and Lemma \ref{lem10}). Unfortinately
this hypothesis was omitted from the statement of theorems in
\cite{ellisJPAA} and \cite{brownellis} and some subsequent
publications. This seems to have caused  confusion. However, in
all calculations based on the theorems the connectivity condition
was met. A detailed erratum \cite{erratum}
 has been available on the first author's home page since 2002.
\end{section}

\begin{section}{Acknowledgments}
We are grateful to R. Brown for his comments on an earlier draft of this paper. We are also grateful to the Max Planck Institute f\"ur Mathematik, Bonn for
inviting both authors to a symposium in honour of H.-J. Baues; work on the paper began during this symposium.

\end{section}

\end{document}